\documentclass{article}[12pt]
\usepackage{amssymb}
\usepackage{epsfig}
\usepackage{amsmath}
\usepackage{amsfonts}
\usepackage{pgf,tikz}
\usepackage{enumerate}
\usepackage[colorlinks=true, linkcolor=black, citecolor=black, urlcolor=black]{hyperref}

\newcommand{\R}{{{\Bbb R}}}

\newcommand{\N}{{{\Bbb N}}}

\newcommand{\T}{{{\mathbb{T}}}}

\newtheorem{theorem}{\sc Theorem}[section]
\newtheorem{proposition}[theorem]{\sc Proposition}
\newtheorem{lemma}[theorem]{\sc Lemma}
\newtheorem{definition}[theorem]{\sc Definition}
\newtheorem{remark}[theorem]{\sc Remark}
\newtheorem{corollary}[theorem]{\sc Corollary}
\newtheorem{example}[theorem]{\sc Example}

\def\qed{\hbox to 0pt{}\hfill$\rlap{$\sqcap$}\sqcup$\medbreak}

\DeclareMathOperator{\co}{co}

\title{A version of Krasnoselskii's compression--expansion fixed point theorem in cones for discontinuous operators with applications}

\author{ Rub\'en Figueroa$^\dag$, Rodrigo L\'opez Pouso$^\dag$
and Jorge Rodr\'{\i}guez-L\'opez}

\date{}
\begin{document}
 \maketitle

\begin{center}  {\small $\dag$ Departamento de Estat\'{\i}stica, An\'alise Matem\'atica e Optimizaci\'on, Universidade de Santiago de Compostela, \\ 15782, Facultade de Matem\'aticas, Campus Vida, Santiago, Spain.}
\end{center}

\medbreak

 \abstract{We introduce a new fixed point theorem of Krasnoselskii type for discontinuous operators. As an application we use it to study the existence of positive solutions of a second--order differential problem with separated boundary conditions and discontinuous nonlinearities.}

\medbreak


\medbreak


\section{Introduction}

A classical problem \cite{Lan,LanWebb,webb_Pos_solut} is that of the existence of positive solutions for the differential equation
\begin{equation}\label{eq_dif_2ord}
u''(t)+g(t)f(u(t))=0, \quad (0<t<1),
\end{equation}
along with suitable boundary conditions (BCs).

This problem arises in the study of radial solutions in $\mathbb{R}^{n}, \ n\geq 2$ for the partial differential equation (PDE)
\[\Delta v+h(\left\|x\right\|)f(v)=0, \quad x\in\mathbb{R}^{n}, \ \left\|x\right\|\in\left[R_{1},R_{2}\right],\]
with the appropriate boundary conditions, see \cite{ErbeHuWang,Lan,LanWebb}.

Recently, in the paper \cite{InfantePietramala2015}, the authors study the existence of non trivial radial solutions for a system of PDEs of the previous type. First, they turn the former problem into a system of ordinary differential equations similar to (\ref{eq_dif_2ord}).

The main novelty in this paper is that we will let $f$ to be discontinuous.

The classical compression--expansion fixed point theorem of Krasnoselskii (see \cite{amann} or \cite{ZeidlerI}) is a well--known tool of nonlinear analysis and it has proved very useful to deduce existence of solutions for nonlinear problems. Here we prove a generalization of that theorem which allows discontinuous operators. The idea is similar to that employed in \cite{figinf,pouso}, where Schauder's fixed point theorem was extended. Then we return to problem (\ref{eq_dif_2ord}) along with Sturm--Liouville BCs and we use our extension of Krasnoselskii's theorem to get a result about existence of positive solutions when $f$ is not necessarily continuous.

\section{Krasnoselskii's fixed point theorem for discontinuous operators}

In the sequel we need the following definitions. A closed and convex subset $K$ of a Banach space $(X,\| \cdot \|)$ is a cone if it satisfies the following conditions:
\begin{enumerate}
	\item[$(i)$] if $x\in K$, then $\lambda x\in K$ for all $\lambda\geq 0$;
	\item[$(ii)$] if $x\in K$ and $-x\in K$, then $x=0$. 
\end{enumerate}
A cone $K$ defines the partial order in $X$ given by $x\preceq y$ if and only if $y-x\in K$.

Let $U$ be a relatively open subset of $K$ and let $T:\overline{U}\subset K\rightarrow K$ be an operator, not necessarily continuous. 
\begin{definition}
	\label{def1}
	The closed--convex envelope of an operator $T:\overline{U} \longrightarrow K$ is the multivalued mapping  $\mathbb{T}: \overline{U} \longrightarrow 2^K$ given by
	\begin{equation}\label{TT}
	\mathbb{T}x=\bigcap_{\varepsilon>0}\overline{\rm co} \, T\left(\overline{B}_{\varepsilon}(x)\cap \overline{U}\right) \quad \text{for every $x\in \overline{U},$}
	\end{equation}
	where $\overline{B}_{\varepsilon}(x)$ denotes the closed ball centered at $x$ and radius $\varepsilon$, and $\overline{\rm co}$ means closed convex hull.  
	
	In other words, we say that $y \in {\mathbb T}x$ if for every $\varepsilon >0$ and every $\rho>0$ there exist $m \in \mathbb N$ and a finite family of vectors $x_i \in \overline{B}_{\varepsilon}(x) \cap \overline{U}$ and coefficients $\lambda_i \in [0,1]$ ($i=1,2,\dots,m$) such that $\sum \lambda_i=1$ and
	$$\left\|y-\sum_{i=1}^m \lambda_i \, Tx_i\right\|< \rho.$$
\end{definition}

The previous definition was formulated for open subsets of a cone, but it works for arbitrary nonempty subsets of
a Banach space (see \cite{pouso}). 

Closed--convex envelopes (cc--envelopes, for short)  need not be upper semicontinuous (usc, for short), see \cite[Example 1.2]{deim}, unless some additional assumptions are imposed on $T$.

\begin{proposition}
	\label{pro1}
	Let $\mathbb T$ be the cc--envelope of an operator $T: \overline{U} \longrightarrow K$. The following properties are satisfied:
	
	\begin{enumerate}
		\item If $T$ maps bounded sets into relatively compact sets, then $\mathbb T$ assumes compact values and it is usc;
		\item If $T \, \overline{U}$ is relatively compact, then $\mathbb T \, \overline{U}$ is relatively compact.
	\end{enumerate}
\end{proposition}

\noindent
{\bf Proof.} Let $x \in \overline{U}$ be fixed and let us prove that $\mathbb T x$ is compact. We know that $\mathbb T x$ is closed, so it suffices to show that it is contained in a compact set. To do so, we note that
$$\mathbb T x=\bigcap_{\varepsilon>0}\overline{\rm co} \, T\left(\overline{B}_{\varepsilon}(x)\cap \overline{U}\right) \subset  \overline{\rm co} \,  T\left(\overline{B}_{1}(x)\cap \overline{U}\right) \subset \overline{\rm co} \, \overline{T\left(\overline{B}_{1}(x)\cap \overline{U}\right)},$$
and $\overline{\rm co} \, \overline{T\left(\overline{B}_{1}(x)\cap \overline{U}\right)}$ is compact because it is the closed convex hull of a compact subset of a Banach space; see \cite[Theorem 5.35]{ali}. Hence $\mathbb T x$ is compact for every $x \in \overline{U}$, and this property allows us to check that $\mathbb T$ is usc by means of sequences, see \cite[Theorem 17.20]{ali}: let $x_n \to x$ in $\overline{U}$ and let $y_n \in \mathbb T x_n$ for all $n \in \mathbb N$ be such that $y_n \to y$; we have to prove that $y \in \mathbb T x$. Let $\varepsilon>0$ be fixed and take $N \in \mathbb N$ such that $\overline{B}_{\varepsilon}(x_n) \subset \overline{B}_{2 \varepsilon}(x)$ for all $n \ge N$. Then we have $y_n \in \overline{\rm co} \, T(\overline{B}_{\varepsilon}(x_n) \cap \overline{U}) \subset  \overline{\rm co} \, T (\overline{B}_{2\varepsilon}(x ) \cap \overline{U})$ for all $n \ge N$, which implies that $y \in \overline{\rm co} \, T (\overline{B}_{2\varepsilon}(x ) \cap \overline{U})$. Since $\varepsilon>0$ was arbitrary, we conclude that $y \in \mathbb T x$.

Arguments are similar for the second part of the proposition. For every $x \in \overline{U}$ and $\varepsilon>0$ we have
$$\overline{\rm co} \, T(\overline{B}_{\varepsilon}(x) \cap \overline{U}) \subset \overline{\rm co} \,  \overline{T \, \overline{U}}  ,$$
and therefore $\mathbb T x \subset \overline{\rm co}   \, \overline{T \, \overline{U}}  $ for all $x \in \overline{U}$. Hence, $\overline{\mathbb{T} \, \overline{U}}$ is compact because it is a closed subset of the compact set $\overline{\rm co}  \,  \overline{T \, \overline{U}} $.
\qed

Now we recall the fixed point theorem mentioned above (see \cite[Theorem 13.D]{ZeidlerI}).

\begin{theorem}[Krasnoselskii]
Let $r_{i}\leq R$ ($i=1,2$) be positive numbers with $r_{1}\neq r_{2}$ and let $T:\overline{B}(0,R)\cap K\rightarrow K$ be a compact mapping. Suppose that
\begin{enumerate}[$(a)$]
	\item $Tx\not\succeq x$ for all $x\in K$ with $\left\|x\right\|=r_{1}$,
	\item $Tx\not\preceq x$ for all $x\in K$ with $\left\|x\right\|=r_{2}$.
\end{enumerate}
Then $T$ has at least a fixed point $x \in K$ such that 
\[\min\left\{r_{1},r_{2}\right\}<\left\|x\right\|<\max\left\{r_{1},r_{2}\right\}.\]
\end{theorem}

In this section we introduce a generalization of the previous theorem which is based on the following idea: given a possibly discontinuous operator $T$, we build its cc--envelope $\T$ and we prove that it has fixed points by means of the version of Krasnoselskii fixed point theorem for multivalued mappings given by Fitzpatrick--Petryshyn \cite{fp_Index}. Then we impose suitable conditions on $T$ which, roughly speaking, guarantee that fixed points of $\T$ are fixed point of $T$ too.

For completeness, we recall \cite[Theorem 3.2]{fp_Index}.

\begin{theorem}\label{th_FP}
Let $X$ be a Fr\'echet space with a cone $K\subset X$. Let $d$ be a metric on $X$ and let $r_{1},r_{2}\in(0,\infty)$, $r=\max\left\{r_{1},r_{2}\right\}$ and $F:\overline{B}(0,r)\cap K\longrightarrow 2^K$ u.s.c. and condensing. Suppose there exists a continuous seminorm $p$ such that $(I-F)\left(\overline{B}(0,r_{1})\cap K\right)$ is $p$-bounded. Moreover, suppose that $F$ satisfies:
\begin{enumerate}
	\item there is some $w\in K$ with $p(w)\neq 0$ and such that $x\not\in F(x)+tw$ for any $t>0$ and $x\in\partial_{K}{B(0,r_1)}$;
	\item $\lambda x\not\in F(x)$ for any $\lambda>1$ and $x\in\partial_{K}{B(0,r_{2})}$.
\end{enumerate}	
Then $F$ has a fixed point $x_0$ with $\min\left\{r_{1},r_{2}\right\}\leq d(x_0,0)\leq\max\left\{r_1,r_2\right\}$.
\end{theorem}	
 	
We are already in a position to introduce and prove the main results in this section, namely, two extensions of Krasnoselskii fixed point theorem for discontinuous operators.

\begin{theorem}\label{th1_kras}
	Let $r_{i}\leq R$ ($i=1,2$) with $r_{1}\neq r_{2}$ positive numbers and $T:\overline{B}(0,R)\cap K\rightarrow K$ a mapping such that $T\left(\overline{B}(0,R)\cap K\right)$ is relatively compact and  
	\begin{equation}\label{cond}
	\left\{x\right\}\cap\T x\subset\left\{Tx\right\} \quad \text{ for all } x\in\overline{B}(0,R)\cap K,
	\end{equation}
	where $\T$ is the cc--envelope of $T$ as defined in (\ref{TT}).
	
	Suppose that
	\begin{enumerate}[$(a)$]
		\item\label{item1_kras1} $\lambda x\not\in\T x$ for all $x\in K$ with $\left\|x\right\|=r_{1}$ and all $\lambda\geq 1$, 
		\item\label{item2_kras1} there exists $w\in K$ with $\left\|w\right\|\neq 0$ such that $x\not\in \T x+\lambda w$ for all $\lambda\geq 0$ and all $x\in K$ with $\left\|x\right\|=r_{2}$.  
	\end{enumerate}
	Then $T$ has at least a fixed point $x \in K$ such that 
	\[\min\left\{r_{1},r_{2}\right\}<\left\|x\right\|<\max\left\{r_{1},r_{2}\right\}.\]
\end{theorem}

\noindent
{\bf Proof.} Notice that the multivalued mapping $\T$ fulfills all the conditions in \hyperref[th_FP]{Theorem \ref{th_FP}}, so there exists a point $x$ such that $x\in\T x$ and \[\min\left\{r_{1},r_{2}\right\}<\left\|x\right\|<\max\left\{r_{1},r_{2}\right\}.\]
Moreover we deduce from (\ref{cond}) that $x=T\,x$ because $\left\{x\right\}\cap\T x=\left\{x\right\}$.
\qed

A second result leans on compression--expansion type conditions.

\begin{theorem}
\label{th2_kras}
	Let $r_{i}\leq R$ ($i=1,2$) with $r_{1}\neq r_{2}$ positive numbers and $T:\overline{B}(0,R)\cap K\rightarrow K$ a mapping such that $T\left(\overline{B}(0,R)\cap K\right)$ is relatively compact and fulfills condition (\ref{cond}).
	
	Let $\T$ be the cc--envelope of $T$ and suppose that
	\begin{enumerate}[$(i)$]
		\item\label{item1_kras2} $y\not\succeq x$ for all $y\in\T\,x$ and all $x\in K$ with $\left\|x\right\|=r_{1}$,
		\item\label{item2_kras2} $y\not\preceq x$ for all $y\in\T\,x$ and all $x\in K$ with $\left\|x\right\|=r_{2}$.
	\end{enumerate}
	Then $T$ has at least a fixed point $x \in K$ such that 
	\[\min\left\{r_{1},r_{2}\right\}<\left\|x\right\|<\max\left\{r_{1},r_{2}\right\}.\]
\end{theorem}

\noindent
{\bf Proof.} It suffices to show that all the conditions in \hyperref[th1_kras]{Theorem \ref*{th1_kras}} are satisfied. First, we show that condition $(\ref{item1_kras2})$ implies condition $(\ref{item1_kras1})$ in Theorem \ref{th1_kras}. Let $x \in K$ be such that  $\left\|x\right\|=r_{1}$ and let $\lambda \ge 1$; we have to prove that $\lambda x \not \in \mathbb T x$. Reasoning by contradiction, we assume that $y=\lambda x \in \mathbb T x$. Then we have
$$y-x=(\lambda-1)x \in K \quad \mbox{(because $\lambda-1 \ge 0$),}$$
and this implies that $y \succeq x$, a contradiction with condition $(\ref{item1_kras2})$. 

Now for condition $(\ref{item2_kras1})$ in Theorem \ref{th1_kras}. Once again we use a contradiction argument: we assume that for every $w \in K$ such that $ \|w\|  \neq 0$ we can find $x\in\partial_{K} B(0,r_{2})$ and $\lambda\geq 0$ such that $x\in\T x+\lambda w$, i.e., there exists $y\in\T\,x$ such that $x=y+\lambda w$. Hence, $x-y=\lambda w\in K$, a contradiction with $(\ref{item2_kras2})$.
\qed

\begin{remark}
Condition (\ref{cond}) is weaker than continuity, since if $T$ is continuous then $\T x=\left\{Tx\right\}$, so (\ref{cond}) is trivially satisfied. In addition, it is not difficult to find discontinuous mappings that verify this condition as we show in our next section.
\end{remark}

Neither Theorem \ref{th1_kras} nor \ref{th2_kras} remain true if we replace $\T$ by $T$ in the assumptions, as we show in the following example.

\begin{example}
In $X=\mathbb{R}^{2}$ we consider the cone $K=\mathbb{R}^{2}_{+}=\left\{(x,y)\in\mathbb{R}^{2}:x,y\geq 0 \right\}$.

Let $0<r<R$ and define a mapping $T:K\rightarrow K$ in polar coordinates as
\[T(\rho,\theta)=\left\{\begin{array}{ll} (0,0), & \text{ if } \rho\neq r, \\[0.2cm] (r,\frac{\pi}{2}), & \text{ if } \theta\in\left[\left.0,\frac{\pi}{4}\right)\right., \ \rho=r, \\[0.2cm] (r,0), & \text{ if } \theta\in\left[\frac{\pi}{4},\frac{\pi}{2}\right], \ \rho=r.  \end{array}\right. \]

Note that $\mathbb T x=\{Tx\}=\{(0,0)\}$ for all $x \in K$ such that $\|x\| \neq r$ because $T$ is continuous at those points. For points $x=(r,\theta)$, with $\theta \in [0,\pi/2]$, we have three possibilities: if $\theta \in [0,\pi/4)$, then $\mathbb T x$ is the segment with endpoints $(0,0)$ and $(r,\pi/2)$; if $\theta \in  (\pi/4,\pi/2]$, then $\mathbb Tx$ is the segment with endpoints $(0,0)$ and $(r,0)$; finally, $\mathbb T (r,\pi/4)$ is the triangle with vertices $(0,0)$, $(r,0)$ and $(r,\pi/2)$. Therefore,
\[\left\{x\right\}\cap\T x\subset\left\{Tx\right\} \quad \text{for all } x\in K.\]
Moreover, conditions (\ref{item1_kras2}) and (\ref{item2_kras2}) in Theorem \ref{th2_kras} are satisfied if we replace $\T$ by $T$ (and we take $r_1=R$ and $r_2=r$). However, $T$ has no fixed point in $\overline{B}(0,R)\setminus B(0,r)$.
\end{example}

\section{Application to Sturm--Liouville problems}

We consider the following generalization of equation (\ref{eq_dif_2ord}) with separated BCs:
\begin{equation}\label{eq_SturmLiouville}
\begin{array}{l} u''(t)+g(t)f(t,u(t))=0 \quad (0<t<1), \\ \alpha u(0)-\beta u'(0)=0, \\ \gamma u(1)+\delta u'(1)=0, \end{array}
\end{equation}
where $\alpha,\beta,\gamma,\delta\geq 0$ and $\Gamma:=\gamma\beta+\alpha\gamma+\alpha\delta>0$.

The usual approach to this problem consists in turning it into a fixed point problem with the integral operator
\[Tu(t):=\int_{0}^{1}{G(t,s)g(s)f(s,u(s))\,ds},\]
where $G$ is the Green's function associated to the differential problem.

Motivated by this situation, we study existence of fixed points of Hammerstein integral operators
\begin{equation}\label{eq_Hammerstein}
Tu(t):=\int_{0}^{1}{k(t,s)g(s)f(s,u(s))\,ds},
\end{equation}
defined in a suitable space. Here we consider $\mathcal{C}([0,1])$, endowed with the usual supremum norm $\left\|u\right\|=\max_{t\in[0,1]}{\left|u(t)\right|}$. 

Fixed points of $T$ will be looked for in the cone
\[K=\left\{u\in\mathcal{C}([0,1]): \ u\geq 0, \ \min_{t\in[a,b]}{u(t)}\geq c\left\|u\right\|\right\},\]
where $[a,b]\subset[0,1]$ and $c\in(0,1]$. This cone was introduced by Guo and it was intensively employed in recent years, for example, see \cite{Inf,Lan,webb_Pos_solut}. 

We suppose that the terms of the Hammerstein equation (\ref{eq_Hammerstein}) satisfy the following hypotheses:
\begin{enumerate}
	\item[(H1)] $f:[0,1] \times [0,\infty) \longrightarrow [0,\infty)$ is such that:
	\begin{enumerate}
	\item Compositions $f(\cdot,u(\cdot))$ are measurable whenever $u \in {\cal C}([0,1])$; and
	\item For each $r>0$ there exists $R>0$ such that $f(t,u) \le R$ for a.a. $t \in [0,1]$ and all $u \in [0,r]$.
	\end{enumerate}
	\item[(H2)] $g$ measurable and $g(s)\geq 0$ almost everywhere.
	\item[(H3)] $k:[0,1]\times[0,1]\rightarrow[0,\infty)$ is continuous.
	\item[(H4)] There exists a measurable function $\Phi:[0,1]\rightarrow[0,\infty)$ satisfying 
	\[\Phi g\in L^{1}(0,1) \text{ and } \int_{a}^{b}{\Phi(s)g(s)\,ds}>0,\] and a constant $c\in(0,1]$ such that
	\[\begin{array}{ll} k(t,s)\leq\Phi(s)  & \text{ for all } t,s\in[0,1], \\ c\,\Phi(s)\leq k(t,s) & \text{ for all } t\in[a,b], \ s\in[0,1]. \end{array}\]
\end{enumerate}

\begin{remark}
Conditions $(H1)-(H4)$ are similar to those requested in \cite{Lan} with the exception that we do not require $f$ to be continuous. In addition, our assumptions are more general than those in \cite{LanWebb} or \cite{webb_Pos_solut} where the authors require $g\in L^{1}(0,1)$ and $\Phi\in\mathcal{C}([0,1])$.
\end{remark}

\begin{lemma}\label{lem_apl_kras_T_comp}
If conditions $(H1)-(H4)$ are satisfied, then the operator $T:K\rightarrow K$ introduced in (\ref{eq_Hammerstein}) is well--defined and maps bounded sets into relatively compact sets.
\end{lemma}

\noindent
{\bf Proof.} The operator $T$ maps $K$ into $K$. Indeed, we have
\[\left\|Tu\right\|=\max_{t\in[0,1]}\left\{\int_{0}^{1}{k(t,s)g(s)f(s,u(s))\,ds}\right\}\leq\int_{0}^{1}{\Phi(s)g(s)f(s,u(s))\,ds}.\]
Moreover, \[\min_{t\in[a,b]}\left\{Tu(t)\right\}\geq c\int_{0}^{1}{\Phi(s)g(s)f(s,u(s))\,ds}.\]
Hence, $Tu\in K$ for every $u\in K$.

Now we prove that if $B \subset K$ is an arbitrary nonempty bounded set, then $T B$ is relatively compact. Let $r>0$ such that $u \in B$ implies $0 \le u(t) \le r$ for all $t \in [0,1]$, and let $R>0$ be the constant associated to $r>0$ by condition $(H1) \, (b)$. Given $u\in B$, we have
\[\int_{0}^{1}{k(t,s)g(s)f(s,u(s))\,ds}\leq R \int_{0}^{1}{\Phi(s)g(s)\,ds}<\infty,\]
so $TB$ is uniformly bounded. To see that $TB$ is equicontinuous, it suffices to show that for every $\tau\in[0,1]$ and $t_n\rightarrow\tau$, we have
\begin{equation}
\label{e1}
\lim_{t_{n}\rightarrow\tau}{\int_{0}^{1}{\left|k(t_{n},s)g(s)f(s,u(s))-k(\tau,s)g(s)f(s,u(s))\right|\,ds}}=0 \quad \mbox{uniformly in $u \in B$.}
\end{equation}
To prove it, we note that for every $u \in B$ we have
\begin{equation}
\label{e2}
\left|k(t_{n},s)g(s)f(s,u(s))-k(\tau,s)g(s)f(s,u(s))\right| \le Rg(s) \left|k(t_{n},s)-k(\tau,s)\right|,
\end{equation}
which tends to zero as $n$ tends to infinity for a.a. $s \in [0,1]$ because $k$ is continuous in $[0,1]$. Moreover, 
\[R \, g(s) \left|k(t_{n},s)-k(\tau,s)\right|\leq 2R\Phi(s)g(s) \quad \text{for all } n\in\mathbb{N},\]
and $2 R \Phi \, g \in L^1(0,1)$, by $(H4)$, so the dominated convergence theorem and (\ref{e2}) yield (\ref{e1}).
\qed

Moreover suppose that the discontinuities of $f$ allow the operator $T$ to satisfy the condition
\begin{equation}\label{condTT}
\left\{u\right\}\cap\T\,u\subset\left\{Tu\right\} \quad \text{ for all } u\in K\cap\T\,K,
\end{equation}
where $\T$ is the multivalued mapping associated to $T$ defined in (\ref{TT}). Examples of this type of nonlinearities $f$ can be looked up in \cite{figinf,pouso}.

\begin{lemma}\label{lem1_apl_Hamm_Kras}
	Suppose that condition (\ref{condTT}) holds and that
	\begin{enumerate}
		\item[$(I_{\rho}^{1})$] There exist $\rho>0$ and $\varepsilon>0$ such that $f^{\rho,\varepsilon}<m$, where 
		\[f^{\rho,\varepsilon}:=\sup_{0 \le t \le 1, \, 0 \le u \le \rho+\varepsilon}\left\{\frac{f(t,u)}{\rho}\right\} \quad \text{ and } \quad \frac{1}{m}:=\sup_{t\in[0,1]}{\int_{0}^{1}{k(t,s)g(s)\,ds}}.\] 
	\end{enumerate}
	Then $\lambda u\not\in\T u$ for all $u\in\partial_{K} B(0,\rho)$ and all $\lambda\geq 1$.
\end{lemma}

\noindent
{\bf Proof.} Suppose that there exist $\lambda\geq 1$ and $u\in\partial_{K} B(0,\rho)$ such that $\lambda u=Tv$ for some $v\in\overline{B}_{\varepsilon}(u)\cap K$, i.e.,
\[\lambda u(t)=\int_{0}^{1}{k(t,s)g(s)f(s,v(s))\,ds}.\]
Taking the supremum for $t\in[0,1]$, 
\begin{equation}\label{eq_lem1_apl_kras}
\lambda\rho\leq\sup_{t\in[0,1]}{\int_{0}^{1}{k(t,s)g(s)f(s,v(s))\,ds}}\leq\rho f^{\rho,\varepsilon}\sup_{t\in[0,1]}{\int_{0}^{1}{k(t,s)g(s)\,ds}} \leq\rho f^{\rho,\varepsilon}\frac{1}{m}<\rho,
\end{equation}
a contradiction.

Given $m\in\mathbb{N}$, it is similarly proved that $\lambda u\neq\sum_{i=1}^{m}{\lambda_{i}Tv_{i}}$ for any $v_{i}\in\overline{B}_{\varepsilon}(u)\cap K$ and $\lambda_{i}\in[0,1]$ with $\sum_{i=1}^{m}{\lambda_{i}}=1$. Hence, $\lambda u\not\in\co\left(T\left(\overline{B}_{\varepsilon}(u)\cap K\right)\right)$.

To see $\lambda u\not\in\overline{\co}\left(T\left(\overline{B}_{\varepsilon}(u)\cap K\right)\right)$ we consider two cases: $\lambda=1$ and $\lambda>1$.

If $\lambda=1$, we have $u\not\in\T u$ because $u\neq Tu$ and $\left\{u\right\}\cap\T u\subset\left\{Tu\right\}$.

If $\lambda>1$, we obtain from (\ref{eq_lem1_apl_kras}) that $\lambda \rho\leq\rho$, that in this case suppose a contradiction too.
\qed

In the sequel we denote
\[V_{\rho}=\left\{u\in K:\min_{a\leq t\leq b}{u(t)}<\rho\right\}.\]
In addition, it is trivial to see that $B(0,\rho)\cap K\subset V_{\rho}\subset B(0,\rho/c)\cap K$, and $V_{\rho}$ is a relatively open subset of $K$ (since minimum function is continuous).

\begin{lemma}\label{lem2_apl_Hamm_Kras}
	Suppose that condition (\ref{condTT}) holds and that
	\begin{enumerate}
		\item[$(I_{\rho}^{0})$] There exist $\rho>0$ and $\varepsilon>0$ such that $f_{\rho,\varepsilon}>M(a,b)$, where 
		\[f_{\rho,\varepsilon}:=\inf_{a\leq t\leq b,\, c(\rho-\varepsilon)\leq u\leq\frac{\rho}{c}+\varepsilon}\left\{\frac{f(t,u)}{\rho}\right\} \quad \text{ and } \quad \frac{1}{M(a,b)}:=\inf_{t\in[a,b]}{\int_{a}^{b}{k(t,s)g(s)\,ds}}.\] 
	\end{enumerate}	
	Then $u\not\in\T u+\lambda e$ for all $u\in\partial\, V_{\rho}$, all $\lambda\geq 0$ and $e(t)\equiv 1$.
\end{lemma} 	

\noindent
{\bf Proof.} Suppose there exist $u\in\partial\,V_{\rho}$ and $\lambda\geq 0$ such that $u=Tv+\lambda e$ for some $v\in\overline{B}_{\varepsilon}(u)\cap K$. Then \[u(t)=\int_{0}^{1}{k(t,s)g(s)f(s,v(s))\,ds}+\lambda.\]
Notice that $\left\|v\right\|\leq\left\|u\right\|+\varepsilon\leq\displaystyle\frac{\rho}{c}+\varepsilon$ and $\min_{t\in[a,b]}{v(t)}\geq c\left\|v\right\|\geq c\left(\left\|u\right\|-\varepsilon\right)\geq c(\rho-\varepsilon)$. Therefore, for $t\in[a,b]$
\begin{equation*}
u(t)=\int_{0}^{1}{k(t,s)g(s)f(s,v(s))\,ds}+\lambda\geq\int_{a}^{b}{k(t,s)g(s)f(s,v(s))\,ds}+\lambda\geq\rho f_{\rho,\varepsilon}\int_{a}^{b}{k(t,s)g(s)\,ds}+\lambda.
\end{equation*}
Taking the infimum in $[a,b]$ we have
\[\rho\geq\rho f_{\rho,\varepsilon}\inf_{t\in[a,b]}{\int_{a}^{b}{k(t,s)g(s)\,ds}}+\lambda>\rho+\lambda,\]
a contradiction because $\lambda\geq 0$.

Given $m\in\N$, it is similar to check that $u\neq\sum_{i=1}^{m}{\lambda_{i}Tv_{i}}+\lambda e$ for any $v_{i}\in\overline{B}_{\varepsilon}(u)$ and $\lambda_{i}\in[0,1]$ ($i=1,\ldots,n$) with $\sum_{i=1}^{m}{\lambda_{i}}=1$. Hence, \[u\not\in\co\left(T\left(\overline{B}_{\varepsilon}(u)\cap K\right)\right)+\lambda e.\]

If we consider two cases: $\lambda=0$ and $\lambda>0$, and we work in a similar way than in the previous lemma we obtain that $u\not\in\T u+\lambda e$.
\qed

\begin{theorem}\label{th_apl_kras_Hamm}
	Under the hypothesis $(H1)$-$(H4)$ and (\ref{condTT}), the Hammerstein integral operator (\ref{eq_Hammerstein}) has at least a positive fixed point in $K$ if either of the following conditions hold:
	\begin{enumerate} [$(a)$]
		\item There exist $\rho_{1},\rho_{2}\in(0,\infty)$ with $\rho_{1}/c<\rho_{2}$ such that $(I_{\rho_{1}}^{0})$ and $(I_{\rho_{2}}^{1})$ hold.
		\item There exist $\rho_{1},\rho_{2}\in(0,\infty)$ with $\rho_{1}<\rho_{2}$ such that $(I_{\rho_{1}}^{1})$ and  $(I_{\rho_{2}}^{0})$ hold.
	\end{enumerate} 
\end{theorem}

\noindent
{\bf Proof.} It is an immediately consequence of the generalization of \hyperref[th1_kras]{Krasnoselskii's Theorem \ref*{th1_kras}} together with both lemmas above: \hyperref[lem1_apl_Hamm_Kras]{Lemma \ref*{lem1_apl_Hamm_Kras}} and \hyperref[lem2_apl_Hamm_Kras]{Lemma \ref*{lem2_apl_Hamm_Kras}}.
\qed

\begin{remark}
Multiplicity results can be obtained combining previous conditions (see \cite{Lan}).
\end{remark}

Now we return to the differential BVP (\ref{eq_SturmLiouville}). We will say that $u$ is a solution of that problem if $u\in W^{2,1}([0,1])$ (i.e, if $u\in\mathcal{C}^{1}([0,1])$ and $u'\in {\rm AC}([0,1])$, where ${\rm AC}([0,1])$ denote the absolutely continuous functions space defined in $[0,1]$) and satisfies (\ref{eq_SturmLiouville}). 

The problem (\ref{eq_SturmLiouville}) was widely studied looking for positive solutions \cite{ErbeHuWang,Lan}. However, the novelty here is to let function $f$ be discontinuous. In \cite{ErbeHuWang}, the authors consider the problem with $g(t)f(t,u)=h(t,u(t))$ where $h$ is continuous and they use a norm compression--expansion theorem in order to guarantee the existence of solutions. On the other hand, in \cite{Lan}, Lan considers $f$ autonomous and continuous and weaker conditions about $g$, he even replaces the hypothesis integrable by measurable, but it is necessary that $\int_{0}^{1}{\Phi(s)g(s)\,ds}<\infty$.
Here, as $f$ can be discontinuous, we will require $g\in L^{1}(0,1)$.

We can write the differential problem (\ref{eq_SturmLiouville}) as
\[u(t)=\int_{0}^{1}{G(t,s)g(s)f(s,u(s))\,ds}=:Tu(t),\]
where $G$ is the associated Green function, that in this case \cite{Lan} is given by
\begin{equation}\label{G_SturmLiouville}
G(t,s)=\frac{1}{\Gamma}\left\{\begin{matrix} (\gamma+\delta-\gamma \,t)(\beta+\alpha\, s), & \text{ if } 0\leq s\leq t\leq 1, \\
(\beta+\alpha\, t)(\gamma+\delta-\gamma\, s), & \text{ if } 0\leq t<s\leq 1, \end{matrix}\right. 
\end{equation}
and it is non negative.

As $G(t,s)\leq G(s,s)$ for all $s,t\in [0,1]$, it is possible to choose \[\Phi(s)=G(s,s)=\frac{1}{\Gamma}(\gamma+\delta-\gamma\,s)(\beta+\alpha\,s).\] 
Moreover we can choose $a$, $b$ and $c$ in the following way \cite{Lan}:
\begin{enumerate}[(C1)]
	\item\label{C1} $a,b\in[0,1]$ such that $-\beta/\alpha<a<b<1+\delta/\gamma$, where we consider $\beta/\alpha=\infty$ if $\alpha=0$ and $\delta/\gamma=\infty$ if $\gamma=0$.
	\item $c=\min\left\{(\gamma+\delta-\gamma\, b)/(\gamma+\delta),(\beta+\alpha\, a)/(\alpha+\beta)\right\}$.
\end{enumerate}
These choices guarantee that $c\,\Phi(s)\leq G(t,s)$ for $t\in[a,b]$ and $s\in[0,1]$.

We shall work, as before, in the cone
\[K=\left\{u\in\mathcal{C}[0,1]: \ u\geq 0, \ \min_{t\in[a,b]}{u(t)}\geq c\left\|u\right\|\right\}.\]

We allow $f:[0,1]\times [0,\infty)\rightarrow[0,\infty)$ to have discontinuities over the graphs of the following curves.

\begin{definition}
	We say that $\gamma:[r,s]\subset I=[0,1]\rightarrow[0,\infty)$, $\gamma\in W^{2,1}([r,s])$, is an admissible discontinuity
	curve for the differential equation $u''=-g(t)f(t,u)$ if one of the following conditions holds:
	\begin{enumerate}[(a)]
		\item $\gamma''(t)=-g(t)f(t,\gamma(t))$ for a.e. $t\in[r,s]$ (then we say $\gamma$ is viable for the differential equation),
		\item There exist $\varepsilon>0$ and $\psi\in L^{1}(r,s), \ \psi(t)>0$ for a.e. $t\in[r,s]$ such that either
		\begin{equation}\label{eq_curva_inviable1}
		\gamma''(t)+\psi(t)<-g(t)f(t,y) \text{ for a.e. } t\in I \text{ and all } y\in\left[\gamma(t)-\varepsilon,\gamma(t)+\varepsilon\right],
		\end{equation}
		or
		\begin{equation}\label{eq_curva_inviable2}
		\gamma''(t)-\psi(t)>-g(t)f(t,y) \text{ for a.e. } t\in I \text{ and all } y\in\left[\gamma(t)-\varepsilon,\gamma(t)+\varepsilon\right].
		\end{equation}
		In this case we say that $\gamma$ is inviable.
	\end{enumerate}	
\end{definition}

Working with admissible discontinuity curves involves some technicalities gathered in the next lemma and its subsequent corollaries whose proofs will be omitted because they can be found in \cite{pouso}.

\begin{lemma}[{\cite[Lemma 4.1]{pouso}}]
	\label{le1}
	Let $a,b \in \R$, $a<b$, and let $g, h \in L^1(a,b)$, $g \ge 0$ a.e., and $h>0$ a.e. in $(a,b)$.
	
	For every measurable set $J \subset (a,b)$ with $m(J)>0$ there is a measurable set $J_0 \subset J$ with $m(J \setminus J_0)=0$ such that for every $\tau_0 \in J_0$ we have
	\begin{equation*}
	\label{lim}
	\lim_{t \to \tau_0^+}\dfrac{\int_{[\tau_0,t]\setminus J}g(s) \, ds}{\int_{\tau_0}^{t}{h(s) \, ds}}=0=\lim_{t \to \tau_0^-}\dfrac{\int_{[t,\tau_0]\setminus J}g(s) \, ds}{\int_{t}^{\tau_0}{h(s) \, ds}}.
	\end{equation*}
\end{lemma}

\begin{corollary}[{\cite[Corollary 4.2]{pouso}}]
	\label{cole}
	Let $a,b \in \R$, $a<b$, and let $h \in L^1(a,b)$ be such that $h>0$ a.e. in $(a,b)$.
	
	For every measurable set $J \subset (a,b)$ with $m(J)>0$ there is a measurable set $J_0 \subset J$ with $m(J \setminus J_0)=0$ such that for all $\tau_0 \in J_0$ we have
	\begin{equation*}
	\label{lim2}
	\lim_{t \to \tau_0^+}\dfrac{\int_{[\tau_0,t] \cap J}h(s) \, ds}{\int_{\tau_0}^{t}{h(s) \, ds}}=1=\lim_{t \to \tau_0^-}\dfrac{\int_{[t,\tau_0]\cap J}h(s) \, ds}{\int_{t}^{\tau_0}{h(s) \, ds}}.
	\end{equation*}
\end{corollary}

\begin{corollary}[{\cite[Corollary 4.3]{pouso}}]
	\label{coabs}
	Let $a,b \in \R$, $a<b$, and let $f, \, f_n:[a,b] \longrightarrow \R$ be absolutely continuous functions on $[a,b]$ ($n \in \N$), such that $f_n \to f$ uniformly on $[a,b]$ and for a measurable set $A \subset [a,b]$ with $m(A)>0$ we have
	$$\lim_{n \to \infty}f_n'(t)=g(t) \quad \mbox{for a.a. $t \in A$.}$$
	
	If there exists $M \in L^1(a,b)$ such that $|f'(t)| \le M(t)$ a.e. in $[a,b]$ and also $|f_n'(t)| \le M(t)$ a.e. in $[a,b]$ ($n \in \N$), then $f'(t)=g(t)$ for a.a. $t \in A$.
\end{corollary}

We shall also need the following result.

\begin{lemma}
\label{lemaQ}
If $M \in L^1(0,1)$, $M \ge 0$ almost everywhere, then the set
\[Q=\left\{u\in\mathcal{C}^{1}([0,1]):\left|u'(t)-u'(s)\right|\leq\int_{s}^{t}{M(r)\,dr} \quad \mbox{whenever $0 \le s\leq t \le 1$} \right\},\]
is closed in $\mathcal{C}([0,1])$ with the maximum norm topology. 

Moreover, if $u_n \in Q$ for all $n \in \mathbb N$ and $u_n\to u$ uniformly in $[0,1]$, then there exists a subsequence $\{u_{n_k}\}$ which tends to $u$ in the ${\cal C}^1$ norm.
\end{lemma}

\noindent
{\bf Proof.} Let $\{u_n\}$ be a sequence of elements of $Q$ which converges uniformly on $[0,1]$ to some function $u \in {\cal C}([0,1])$; we have to show that $u \in Q$ and a subsequence $\{u_{n_k}\}$ tends to $u$ in the ${\cal C}^1$ norm. 

Since each $u_n$ is continuously differentiable, the Mean Value Theorem guarantees the existence of some $t_n \in (0,1)$ such that
$$u_n'(t_n)=u_n(1)-u_n(0).$$
This implies the existence of some $K>0$ such that $|u_n'(t_n)| \le K$ for all $n \in \mathbb N$, because $\{u_n\}$ is uniformly bounded in $[0,1]$. Hence, for every $n \in \mathbb N$ and every $t \in [0,1]$, we have
$$|u_n'(t)| \le |u_n'(t)-u_n'(t_n)|+|u_n'(t_n)| \le \int_0^1{M(s) \, ds}+K,$$
so $\{u_n\}$ is bounded in the ${\cal C}^1$ norm. Moreover, the definition of $Q$ implies that the sequence $\{u_n'\}$ is equicontinuous in $[0,1]$, so the Ascoli--Arzel\'a Theorem ensures that some subsequence of $\{u_n\}$, say $\{u_{n_k}\}$, which converges in the ${\cal C}^1$ norm to some $v \in {\cal C}^1([0,1])$. As a result, $u=v$, so $u$ is continuously differentiable in $[0,1]$ and $\{u_{n_k}\}$ tends to $u$ in the ${\cal C}^1$ norm. In particular, $\{u_{n_k}'\}$ tends to $u'$ uniformly in $[0,1]$.

Moreover, for $s,t \in [0,1]$, $s \le t$, and all $k \in \mathbb N$, we have
$$|u_{n_k}'(t)-u_{n_k}'(s)| \le \int_s^t{M(r) \, dr},$$
and going to the limit as $k$ tends to infinity we deduce that $|u'(t)-u'(s)| \le \int_s^t{M(r) \, dr}$.\qed

We are now ready for the proof of the main result in this section.

\begin{theorem}
	Suppose that $f$ and $g$ satisfy the following hypothesis:
	\begin{enumerate}[i.]
		\item $f:[0,1] \times [0,\infty) \longrightarrow [0,\infty)$ is such that:
		\begin{enumerate}
			\item Compositions $f(\cdot,u(\cdot))$ are measurable whenever $u \in {\cal C}([0,1])$; and
			\item For each $r>0$ there exists $R>0$ such that $f(t,u) \le R$ for a.a. $t \in [0,1]$ and all $u \in [0,r]$.
		\end{enumerate}
		\item There exist admissible discontinuity curves $\gamma_{n}:I_{n}=[a_{n},b_{n}]\rightarrow[0,\infty)$, $n\in\N$, such that for a.a. $t \in I$ the function $u \mapsto f(t,u)$ is continuous on $[0,\infty)\setminus\bigcup_{\left\{n:t\in I_n\right\}}\left\{\gamma_{n}(t)\right\}$. 
		\item $g\in L^{1}(0,1)$ and $g(s)\geq 0$ almost everywhere with $\int_{a}^{b}{g(s)\,ds}>0$, where $a$ and $b$ are given in (C\ref{C1}).
	\end{enumerate}
	Moreover, assume that one of the following conditions hold:
	\begin{enumerate} [$(a)$]
		\item There exist $\rho_{1},\rho_{2}\in(0,\infty)$ with $\rho_{1}/c<\rho_{2}$ such that $(I_{\rho_{1}}^{0})$ and $(I_{\rho_{2}}^{1})$ hold.
		\item There exist $\rho_{1},\rho_{2}\in(0,\infty)$ with $\rho_{1}<\rho_{2}$ such that $(I_{\rho_{1}}^{1})$ and  $(I_{\rho_{2}}^{0})$ hold.
	\end{enumerate}
	Then the differential problem with separated BCs (\ref{eq_SturmLiouville}) has at least a positive solution $u\in W^{2,1}([0,1])$.
\end{theorem}

\noindent
{\bf Proof.} The operator $T:K\rightarrow K$ given by
\[Tu(t)=\int_{0}^{1}{G(t,s)g(s)f(s,u(s))\,ds}\]
is well defined and it maps bounded sets into relatively compact ones, as consequence of \hyperref[lem_apl_kras_T_comp]{Lemma \ref*{lem_apl_kras_T_comp}}. In addition, as $G$ is the Green function associated to a second--order homogeneous differential problem, $Tu\in W^{2,1}([0,1])$ for all $u\in K$. On the other hand, given $u\in B(0,\rho_{2}/c)\cap K=K_{2}$, we have $g(t)f(t,u(t))\in L^{1}[0,1]$, and there exists $M(t)\in L^{1}[0,1]$ such that
\begin{equation}\label{CM} 
h(t,u):=g(t)f(t,u)\leq M(t) \text{ for a.e. } t\in[0,1] \text{ and all } u\in K_{2}.
\end{equation}
We consider the set 
\begin{equation}
\label{eQ}
Q=\left\{u\in\mathcal{C}^{1}([0,1]):\left|u'(t)-u'(s)\right|\leq\int_{s}^{t}{M(r)\,dr} \quad (s\leq t) \right\},
\end{equation}
which is  closed in $(\mathcal{C}([0,1]),\left\|\cdot\right\|_{\infty})$ by virtue of Lemma \ref{lemaQ}.

Hence, since $T\,K_{2}\subset Q$ and $Q$ is a closed and convex subset of $\mathcal{C}([0,1])$, we have $\T\,K_{2}\subset Q$. 

Now we will prove that
\begin{equation}\label{condTT_proof}
\left\{u\right\}\cap\T u\subset\left\{Tu\right\} \quad \text{ for all } u\in K_{2}\cap\T K_{2}.
\end{equation}

To do so, we fix an arbitrary function $u \in K_{2}\cap Q$ and we consider three different cases.


\noindent
{\it Case 1 -- $m(\{t \in I_n \, : \, u(t)=\gamma_n(t)\})=0$ for all $n \in \N$}. Let us prove that then $T$ is continuous at $u$. 


The assumption implies that for a.a. $t \in I$ the mapping $h(t,\cdot)$ is continuous at $u(t)$. Hence if $u_k \to u$ in $K_{2}\cap Q$ then
$$h(t,u_k(t)) \to h(t,u(t)) \quad \mbox{for a.a. $t \in I$,}$$
which, along with (\ref{CM}), yield $Tu_k \to T u$ in ${\cal C}(I)$.

\noindent
{\it Case 2 -- $m(\{t \in I_n \, : \, u(t)=\gamma_n(t)\})>0$ for some $n \in \N$ such that $\gamma_n$ is inviable.} In this case we can prove that $u \not \in {\mathbb T}u$.


First, we fix some notation. Let us assume that for some $n \in \N$ we have $m(\{t \in I_n \, : \, u(t)=\gamma_n(t)\})>0$ and there exist $\varepsilon>0$
and $\psi \in L^1(I_n)$, $\psi(t)>0$ for a.a. $t \in I_n$, such that (\ref{eq_curva_inviable2}) holds with $\gamma$ replaced by $\gamma_n$. (The proof is similar if we assume (\ref{eq_curva_inviable1}) instead of (\ref{eq_curva_inviable2}), so we omit it.)

We denote $J=\{t \in I_n \, : \, u(t)=\gamma_n(t)\}$, and we deduce from \hyperref[le1]{Lemma \ref{le1}} that there is a measurable set $J_0 \subset J$ with $m(J_0)=m(J)>0$ such that for all $\tau_0 \in J_0$ we have
\begin{equation}
\label{limi}
\lim_{t \to \tau_0^+}\dfrac{2\int_{[\tau_0,t]\setminus J}M(s) \, ds}{(1/4)\int_{\tau_0}^{t}{\psi(s) \, ds}}=0=\lim_{t \to \tau_0^-}\dfrac{2\int_{[t,\tau_0]\setminus J}M(s) \, ds}{(1/4)\int_{t}^{\tau_0}{\psi(s) \, ds}}.
\end{equation}
By \hyperref[cole]{Corollary \ref{cole}} there exists $J_1 \subset J_0$ with $m(J_0 \setminus J_1)=0$ such that for all $\tau_0 \in J_1$ we have
\begin{equation}
\label{limi2}
\lim_{t \to \tau_0^+}\dfrac{\int_{[\tau_0,t] \cap J_0}\psi(s) \, ds}{\int_{\tau_0}^{t}{\psi(s) \, ds}}=1=\lim_{t \to \tau_0^-}\dfrac{\int_{[t,\tau_0]\cap J_0}\psi(s) \, ds}{\int_{t}^{\tau_0}{\psi(s) \, ds}}.
\end{equation}

Let us now fix a point $\tau_0 \in J_1$. From (\ref{limi}) and (\ref{limi2}) we deduce that there exist $t_-<\tilde{t}_{-}<\tau_0$ and $t_+>\tilde{t}_{+}>\tau_0$, $t_{\pm}$ sufficiently close to $\tau_0$ so that the following inequalities are satisfied for all $t\in[\tilde{t}_{+},t_{+}]$:
\begin{eqnarray}
\label{in+1}
2\int_{[\tau_0,t]\setminus J}M(s) \, ds < \dfrac{1}{4}\int_{\tau_0}^{t}{\psi(s) \, ds}, \\
\label{in+2}
\int_{[\tau_0,t] \cap J}\psi(s) \, ds \ge \int_{[\tau_0,t] \cap J_0}\psi(s) \, ds>\dfrac{1}{2}
\int_{\tau_0}^{t}\psi(s) \, ds, 
\end{eqnarray}
and for all $t\in[t_{-},\tilde{t}_{-}]$:
\begin{eqnarray}
\label{in-1}
2\int_{[t,\tau_0]\setminus J}M(s) \, ds < \dfrac{1}{4}\int_{t}^{\tau_0}{\psi(s) \, ds}, \\
\label{in-2}
\int_{[t,\tau_0] \cap J}\psi(s) \, ds>\dfrac{1}{2}
\int_{t}^{\tau_0}\psi(s) \, ds.
\end{eqnarray}

Finally, we define a positive number
\begin{equation}
\label{ro}
\tilde\rho=\min \left\{ \dfrac{1}{4}\int_{\tilde{t}_-}^{\tau_0}{\psi(s) \, ds}, \dfrac{1}{4}\int_{\tau_0}^{\tilde{t}_+}{\psi(s) \, ds} \right\},
\end{equation}
and we are now in a position to prove that $u \not \in {\mathbb T}u$. It suffices to prove the following claim:

\noindent
{\it Claim -- Let $\varepsilon>0$ be given by our assumptions over $\gamma_n$ and let $\rho=\displaystyle\frac{\tilde{\rho}}{2}\min\left\{\tilde{t}_{-}-t_{-},t_{+}-\tilde{t}_{+}\right\}$ be where $\tilde{\rho}$ is as in (\ref{ro}). For every finite family $u_i \in \overline{B}_{\varepsilon}(x) \cap K$ and $\lambda_i \in [0,1]$ ($i=1,2,\dots,m$), with $\sum \lambda_i=1$, we have $\|u-\sum\lambda_i Tu_i\|\ge \rho$.}


Let $u_i$ and $\lambda_i$ be as in the Claim and, for simplicity, denote $y=\sum\lambda_i Tu_i$. For a.a. $t \in J=\{t \in I_n \, : \, u(t)=\gamma_n(t)\}$ we have
\begin{equation}
\label{etiq}
y''(t)=\sum_{i=1}^m \lambda_i (Tu_i)''(t)=-\sum_{i=1}^m{\lambda_i \, h(t,u_i(t))}.
\end{equation}
On the other hand, for every $i \in \{1,2,\dots,m\}$ and every $t \in J$ we have 
$$|u_i(t)-\gamma_n(t)|=|u_i(t)-u(t)|<\varepsilon,$$
and then the assumptions on $\gamma_n$ ensure that for a.a. $t \in J$ we have
\begin{equation}
\label{et}
y''(t)=-\sum_{i=1}^m{\lambda_i \, h(t,u_i(t))}<\sum_{i=1}^m{\lambda_i \, (\gamma_n''(t)-\psi(t))}=u''(t)-\psi(t).
\end{equation}

Now for $t\in[t_{-},\tilde{t}_{-}]$ we compute
\begin{align*}
y'(\tau_0)-y'(t)&=\int_{t}^{\tau_0}{y''(s) \, ds}=\int_{[t,\tau_0]\cap J}{y''(s) \, ds}+\int_{[t,\tau_0]\setminus J}{y''(s) \, ds}\\
&<\int_{[t,\tau_0]\cap J}{u''(s) \, ds}-\int_{[t,\tau_0]\cap J}{\psi(s) \, ds}\\
& \qquad +\int_{[t,\tau_0]\setminus J}{M(s) \, ds} \qquad (\mbox{by (\ref{et}), (\ref{etiq}) and (\ref{CM})})\\
&=u'(\tau_0)-u'(t)-\int_{[t,\tau_0]\setminus J}{u''(s) \, ds}-\int_{[t,\tau_0]\cap J}{\psi(s) \, ds} +\int_{[t,\tau_0]\setminus J}{M(s) \, ds}\\
&\le u'(\tau_0)-u'(t)-\int_{[t,\tau_0]\cap J}{\psi(s) \, ds}+2\int_{[t,\tau_0]\setminus J}{M(s) \, ds} \\
&<u'(\tau_0)-u'(t)-\dfrac{1}{4}\int_{t}^{\tau_0}{\psi(s) \, ds} \quad (\mbox{by (\ref{in-1}) and (\ref{in-2})}),
\end{align*}
hence $y'(t)-u'(t) \ge \tilde\rho$ provided that $y'(\tau_0) \ge u'(\tau_0)$. Therefore, by integration we obtain
\[y(\tilde{t}_{-})-u(\tilde{t}_{-})=y(t_-)-u(t_-)+\int_{t_-}^{\tilde{t}_{-}}{(y'(t)-u'(t))\,dt}\geq y(t_-)-u(t_-)+\tilde{\rho}(\tilde{t}_{-}-t_{-}),\]
so if $y(t_{-})-u(t_{-})\leq -\rho$, then $\left\|y-u\right\|\geq\rho$. Otherwise, if $y(t_{-})-u(t_{-})>-\rho$, then we have $y(\tilde{t}_{-})-u(\tilde{t}_{-})>\rho$ and thus $\left\|y-u\right\|\geq\rho$ too.

Similar computations in the interval $[\tilde{t}_{+},t_+]$ instead of $[t_-,\tilde{t}_{-}]$ show that if $y'(\tau_0) \le u'(\tau_0)$ then we have $u'(t)-y'(t) \ge \tilde\rho$ for all $t\in[\tilde{t}_{+},t_+]$ and this also implies $\left\|y-u\right\|\geq\rho$. 
The claim is proven.

\noindent
{\it Case 3 -- $m(\{t \in I_n \, : \, u(t)=\gamma_n(t)\})>0$ only for some of those $n \in \N$ such that $\gamma_n$ is viable.} Let us prove that in this case the relation $u \in {\mathbb T}u$ implies $u=Tu$.

Let us consider the subsequence of all viable admissible discontinuity curves in the conditions of Case 3, which we denote again by $\{\gamma_n\}_{n \in {\scriptsize \N}}$ to avoid overloading notation. We have $m(J_n)>0$ for all $n \in \N$, where
$$J_n=\{t \in I_n \, : \, u(t)=\gamma_n(t)\}.$$
For each $n \in \N$ and for a.a. $t \in J_n$ we have 
$$u''(t)=\gamma_n''(t)=-h(t,\gamma_n(t))=-h(t,u(t)),$$
and therefore $u''(t)=-h(t,u(t))$ a.e. in $J=\cup_{n \in \scriptsize{\N}}J_n$.

Now we assume that $u \in {\mathbb T}u$ and we prove that it implies that $u''(t)=-h(t,u(t))$ a.e. in $I \setminus J$, thus showing that $u=Tu$.

Since $u \in {\mathbb T}u$ then for each $k \in \N$ we can guarantee that we can find functions $u_{k,i} \in \overline{B}_{1/k}(u)\cap K_{2}$ and coefficients $\lambda_{k,i} \in [0,1]$ ($i=1,2,\dots,m(k)$) such that $\sum \lambda_{k,i}=1$ and
$$\left\|u-\sum_{i=1}^{m(k)}\lambda_{k,i}Tu_{k,i}\right\|<\dfrac{1}{k}.$$

Let us denote $y_k=\sum_{i=1}^{m(k)}\lambda_{k,i}Tu_{k,i}$, and notice that $y_k \to u$ uniformly in $I$ and $\|u_{k,i}-u\|\le 1/k$ for all $k \in \N$ and all $i \in \{1,2,\dots,m(k)\}$.

For every $k\in\mathbb{N}$ we have $y_k\in Q$ as defined in (\ref{eQ}), and therefore Lemma \ref{lemaQ} guarantees that $u \in Q$ and, up to a subsequence, $y_k \to u$ in the ${\cal C}^1$ topology.

For a.a. $t \in I\setminus J$ we have that $h(t,\cdot)$ is continuous at $u(t)$, so for any $\varepsilon>0$ there is some $k_0=k_0(t) \in \N$ such that for all $k \in \N$, $k \ge k_0$, we have
$$|h(t,u_{k,i}(t))-h(t,u(t))| < \varepsilon \quad \mbox{for all $i \in \{1,2,\dots,m(k)\}$,}$$
and therefore
$$|y_k''(t)+h(t,u(t))|\le \sum_{i=1}^{m(k)}\lambda_{k,i}|h(t,u_{k,i}(t))-h(t,u(t))| <\varepsilon.$$
Hence $y_k''(t) \to -h(t,u(t))$ for a.a. $t \in I \setminus J$, and then Corollary \ref{coabs} guarantees that $u''(t)=-h(t,u(t))$ for a.a. $t \in I \setminus J$. 

Therefore the proof of condition (\ref{condTT_proof}) is over and we conclude by means of \hyperref[th_apl_kras_Hamm]{Theorem \ref{th_apl_kras_Hamm}}.
\qed

\begin{remark}
The differential problem (\ref{eq_SturmLiouville})  contains Dirichlet and Robin problems, so the previous result generalizes the existence results given in \cite{LanWebb}, because here we allow $f$ be discontinuous.
\end{remark}

%
%

\section*{Competing interests}
The authors declare that they have no competing interests.

\section*{Acknowledgements}
Rodrigo L\'opez Pouso was partially supported by
Ministerio de Econom\'{\i}a y Competitividad, Spain, and FEDER, Project
MTM2013-43014-P, and Xunta de Galicia R2014/002 and GRC2015/004. Rub\'en Figueroa and Jorge Rodr\'iguez-L\'opez were partially supported by Xunta de Galicia, project EM2014/032.

\end{document}